\documentclass[11pt]{article}

\usepackage{color,latexsym,amsmath,amssymb,path,enumerate,url}

\newcommand{\ignore}[1]{}

\newtheorem{theorem}{Theorem}[section]

\newtheorem{lemma}[theorem]{Lemma}

\newtheorem{proposition}[theorem]{Proposition}

\newcommand{\Proof}[1]
{\noindent\emph{Proof #1.}~}

\newsavebox{\smallProofsym}                     % smallproofsym.tex

\savebox{\smallProofsym}                        %
        {%                                      %
        \begin{picture}(7,7)                    %
        \put(0,0){\framebox(6,6){}}             %
        \put(0,2){\framebox(4,4){}}             %
        \end{picture}                           %
        }                                       %

\newcommand{\smalleop}[1]
        {
        \mbox{} \hfill #1~~\usebox{\smallProofsym}\!\!\!\!\!\!\
        }

\newenvironment{theProof}[1]
        {
        \Proof{#1}}{\smalleop{}
        \medskip
        }

\usepackage{graphicx,psfrag}

\usepackage[rflt]{floatflt}

\newcommand{\placefig}[2]
        {\includegraphics[width=#2]{#1.eps}}

\usepackage{dsfont}

\newcommand{\NN}{\ensuremath{\mathbb N}}

\newcommand{\RR}{\ensuremath{\mathbb R}}

\newcommand{\pts}{\mathcal P}

\begin{document}

\pagenumbering{arabic}

\title{On lattices, distinct distances, and the Elekes-Sharir framework\thanks{%
Work by Javier Cilleruelo has been supported by grants MTM 2011-22851
of MICINN and ICMAT Severo Ochoa project SEV-2011-0087.
Work by Adam Sheffer and Micha Sharir has been supported
by Grant 338/09 from the Israel Science Fund,
by the Israeli Centers of Research Excellence (I-CORE)
program (Center No.~4/11), and
by the Hermann Minkowski-MINERVA Center for Geometry
at Tel Aviv University. }}

\author{
Javier Cilleruelo\thanks{%
Instituto de Ciencias Matematicas (CSIC-UAM-UC3M-UCM), and
Departamento de Matematicas, Universidad Aut\'onoma de Madrid, 28049 Madrid, Spain.
{\sl franciscojavier.cilleruelo@uam.es} }
\and
Micha Sharir\thanks{%
School of Computer Science, Tel Aviv University,
Tel Aviv 69978, Israel.
{\sl michas@tau.ac.il} }
\and
Adam Sheffer\thanks{%
School of Computer Science, Tel Aviv University,
Tel Aviv 69978, Israel.
{\sl sheffera@tau.ac.il}
}}

%\date{}

\maketitle

\begin{abstract}
In this note we consider distinct distances determined by points in an
integer lattice. We first consider Erd{\H o}s's lower bound for the square
lattice, recast in the setup of the so-called Elekes-Sharir framework
\cite{ES11,GK11}, and show that, without a major change, this framework
\emph{cannot} lead to Erd{\H o}s's conjectured lower bound. This shows that the
upper bound of Guth and Katz \cite{GK11} for the related 3-dimensional
line-intersection problem is tight for this instance. The gap between this bound
and the actual bound of Erd{\H o}s arises from an application of the
Cauchy-Schwarz inequality (which is an integral part of the Elekes-Sharir
framework). Our analysis relies on two number-theoretic results by Ramanujan.

We also consider distinct distances in rectangular lattices of the form
$\{(i,j) \mid 0\le i\le n^{1-\alpha},\ 0\le j\le n^{\alpha}\}$,
for some $0<\alpha<1/2$, and show that the number of distinct distances in such
a lattice is $\Theta(n)$. In a sense, our proof ``bypasses" a deep conjecture in number theory, posed by Cilleruelo and Granville \cite{CG07}. A positive resolution of this conjecture would also have implied our bound.
\end{abstract}

\section{On the limitations of the Elekes-Sharir framework}\label{sec:limitation}
Given a set $\pts$ of $n$ points in $\RR^2$, let $D(\pts)$ denote the number
of distinct distances that are determined by pairs of points from $\pts$. Let $D(n) = \min_{|\pts|=n}D(\pts)$;
that is, $D(n)$ is the minimum number of distinct distances that any set of $n$ points in $\RR^2$ must always determine.
In his celebrated 1946 paper \cite{erd46}, Erd\H os derived the bound $D(n) = O(n/\sqrt{\log n})$ by considering a $\sqrt{n}\times \sqrt{n}$ integer lattice.
Recently, after 65 years and a series of increasingly larger lower bounds\footnote{For a comprehensive list of the previous bounds, see \cite{GIS11} and
\url{http://www.cs.umd.edu/~gasarch/erdos_dist/erdos_dist.html} (version of May 2013).}, Guth and Katz \cite{GK11} provided an almost matching lower bound
$D(n) = \Omega(n/\log n)$. For this, Guth and Katz simplified the Elekes-Sharir framework \cite{ES11},
developed by Elekes, and slightly refined by Sharir, for tackling the distinct distances problem,
to make it reduce this latter problem to a problem about line intersections in $\RR^3$.
To solve this problem, Guth and Katz developed several novel techniques, relying on tools from algebraic geometry and 19th century analytic geometry.

In this note, we examine the gap of $\Theta(\sqrt{\log n})$ between Erd\H os's upper bound and Guth and Katz's lower
bound for the vertex set of the square lattice considered by
Erd\H os. While it is conceivable that there exists a set of $n$ points that spans $\Theta(n/\log n)$ distinct distances, the common belief is that $D(n) =
\Theta(n/\sqrt{\log n})$. We prove that, even if this common belief
is correct, the Elekes-Sharir framework cannot lead to the actual bound without a major
improvement of the technique (more specifically, without replacing the use of the Cauchy-Schwarz inequality
by some other technique).

We begin by recalling some of the basics of the Elekes-Sharir framework \cite{ES11,GK11}. Consider a set $\pts$ of $n$ points in the plane and put $x =
D(\pts)$. The reduction revolves around the set
\[ Q = \left\{ (a,p,b,q)\in \pts^4 \mid |ap|=|bq|>0 \right\}; \]
The quadruples in $Q$ are ordered, in the sense that $(a,p,b,q)$, $(b,p,a,q)$, $(p,a,q,b)$, and the other possible
permutations are all considered as distinct elements of $Q$.
Also, in a quadruple $(a,p,b,q)\in Q$, the segments $ap$ and $bq$ are allowed to share vertices
or even coincide.

Basically, the reduction is just double counting $|Q|$, and we now present the lower bound.
We put $x=D(\pts)$, and denote the set of (nonzero) distinct distances that are determined by $\pts\times\pts$ as $\delta_1,\ldots,\delta_x$.
Also, for $1 \le i \le x$, we set
\[ E_i = \left\{(p,q)\in \pts^2 \mid |pq| = \delta_i \right\}.\]
As before, we consider $(p,q)$ and $(q,p)$ as two distinct pairs in $E_i$.
Notice that $\sum_{i=1}^x|E_i|=n^2-n$ since every ordered pair of distinct points of $\pts\times\pts$ is contained
in a unique set $E_i$.  By applying the Cauchy-Schwarz inequality, we have

\begin{equation} \label{eq:lowQ}
|Q| = \sum_{i=1}^x |E_i|^2 \ge \frac{1}{x} \left(\sum_{i=1}^x |E_i| \right)^2 = \frac{(n^2-n)^2}{x}.
\end{equation}

Guth and Katz \cite{GK11} derive the upper bound $|Q| = O(n^3\log n)$ which,
combined with the above lower bound, immediately implies $x = \Omega(n/\log n)$. Deriving
this upper bound is considerably more complicated, and we do not discuss it here.
We show that, for the vertex set of the square lattice, Guth and Katz's bound is tight; that is, $|Q|=\Theta(n^3\log n)$.

\begin{figure}[h]\centerline{\placefig{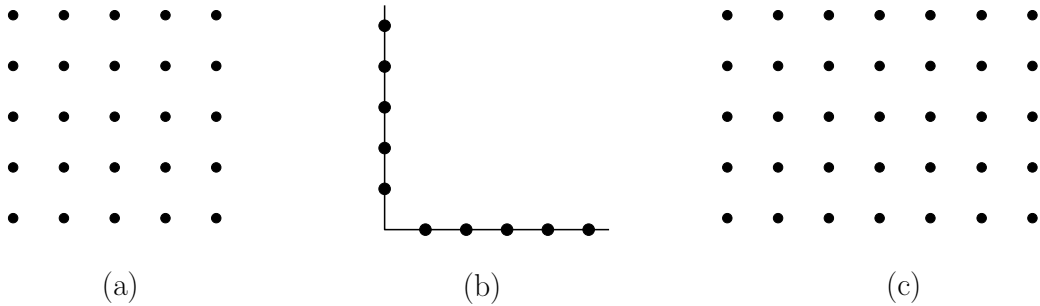}{0.9\textwidth}}
\vspace{-1mm}
\caption{\small \sf (a) A $\sqrt{n}\times\sqrt{n}$ integer lattice.
(b) An L-shaped configuration. (c) A non-square integer lattice.}
\label{fi:config}
\vspace{-2mm}
\end{figure}

We next recall the full details of Erd\H os's lower bound \cite{erd46}.
We take $\pts = \{ (i,j) \mid 1 \le i,j \le \sqrt{n}\}$; that is, $\pts$ is the vertex set of a $\sqrt{n}\times\sqrt{n}$ integer lattice, as depicted in Figure
\ref{fi:config}(a). We require the following theorem from number theory.

\begin{theorem} \label{th:LandRam} {\bf (Landau-Ramanujan \cite{BR95,Land08})}
The number of positive integers smaller than $n$ that are the sum of two squares is $\Theta(n/\sqrt{\log n})$.
\end{theorem}
Every distance determined by a pair from $\pts\times\pts$ is of the form $\sqrt{i^2 + j^2}$ where $0 \le i,j \le \sqrt{n}$. The cases where $i=0$ or $j=0$
contribute negligible amounts to $D(\pts)$ and may thus be ignored. Therefore, Theorem \ref{th:LandRam} implies that $D(\pts) = \Theta(n/\sqrt{\log n})$.
This in turn implies that, for this special set $\pts$, the value of the rightmost expression in (\ref{eq:lowQ}) is $\Theta(n^3\sqrt{\log n})$. Notice that
there is a gap of $\sqrt{\log n}$ between this bound and the upper bound $|Q|=O(n^3\log n)$. To see where this gap comes from, we rely on another result by
Ramanujan.

Recall the expression $\sum_{i=1}^x |E_i|^2$ from (\ref{eq:lowQ}) just before applying the Cauchy-Schwarz inequality. For a positive
integer $k$, we set $r(k) = \left|\{(i,j) \in \NN^2 \mid i^2+j^2=k\}\right|$, and also $\hat{r}(k) = \sum_{i=1}^k r(i)^2$. Consider a distance $1\le \delta
\le \sqrt{n}/10$, say, such that $\delta = \sqrt{i^2 +j^2}$ for some $i,j\in \NN$. Clearly, every point of $\pts$ lies at distance $\delta$ from $r(\delta^2)$ other points of $\pts$.
This implies that
\[\sum_{i:\delta_i<\sqrt{n}/10}|E_i|^2 = \Omega(n^2\hat{r}(n/100)).\]

\begin{theorem} \label{th:Ram} {\bf (Ramanujan \cite{Ram16})}.
$\displaystyle \hat{r}(k) = \Theta\left(k\log k\right)$.
\end{theorem}
By combining Theorem \ref{th:Ram} with the above reasoning, we obtain
\[ \sum_{i=1}^x |E_i|^2 = \Omega(n^2\hat{r}(n/100)) = \Omega(n^3\log n).\]

Thus, by applying the Cauchy-Schwarz inequality, we lose a factor of $\sqrt{\log n}$ which prevents the framework from implying the correct bound
$D(\pts)=\Theta(n/\sqrt{\log n})$. Moreover, since we obtain the lower bound of $|Q|=\Omega(n^3\log n)$, Guth and Katz's upper bound on the number of line
intersections is tight in this case.

\paragraph{Remarks.}(i) The fact that the Cauchy-Schwarz inequality forms a gap in this case indicates that the sizes of the sets $E_i$ are highly
non-uniform. Note that $E_i = \Theta(nr(\delta_i^2))$ (as long as $\delta_i$ is not too close to $\sqrt{n}$).
Since there are only $x=\Theta(n/\sqrt{\log n})$ such sets (by Theorem \ref{th:LandRam}), and since $\sum_{i=1}^x |E_i| = \Theta(n^2)$, the average size of
these sets is $\Theta(n\sqrt{\log n})$.
That is, the average value of the quantities $r(\delta_i^2)$, for $i=1,\ldots,x$, is $\Theta(\sqrt{\log n})$, but the average of their squares is
$\displaystyle \frac{1}{x} \sum_{i=1}^x r(\delta_i^2)^2 = \Theta\left((\log n)^{3/2}\right)$
(by Theorem \ref{th:Ram}). This indicates that the values
$r(\delta_i^2)$ form a rather non-uniform sequence, for which the Cauchy-Schwarz inequality is too weak.
\medskip

\noindent (ii) One might also compare this with Erd\H os's lower bound construction for the \emph{repeated distances problem}
\cite{erd46}, which relied on the property that for
infinitely many values $r(\delta_i^2)$ can be as high as $n^{c/\log\log n}$ (for some universal constant $c$). \medskip

\paragraph{The case of an $L$-shaped configuration.}
We next consider a different configuration, in which the Cauchy-Schwarz inequality yields a surprisingly large gap
of $\Theta(n/\sqrt{\log n})$. We set
$\pts_1 = \{(1,0),(2,0),\ldots,(n,0)\}$, a set of $n$ evenly spaced points on the $x$-axis, and similarly, define $\pts_2 =
\{(0,1),(0,2),\ldots,(0,n)\}$ on the $y$-axis. Figure \ref{fi:config}(b) illustrates the configuration $\pts' = \pts_1\cup\pts_2$. Every distance determined by $\pts'$ is
the square root of a sum of two squares, where each such sum is between 0 and $2n^2$.
Thus, Theorem \ref{th:LandRam} implies $D(\pts') = \Theta(n^2/\sqrt{\log n})$.

We define $Q$ as before and repeat the analysis in (\ref{eq:lowQ}), which implies
the lower bound $|Q| = \Omega(n^4/D(\pts')) = \Omega(n^2\sqrt{\log n})$.
However, the value of the expression from (\ref{eq:lowQ}) before applying the
Cauchy-Schwarz inequality is $\Theta(n^3)$. To see why, let $d_i$, for $i=1,\ldots,n-1$, denote the
number of pairs in $\pts'^2$ that span a distance of $i$ (not to be confused with $\delta_i$).

For $1\le i \le n/2$, we have $d_i =\Theta(n)$ (most of these distances are realized in
$(\pts_1\times \pts_1) \cup (\pts_2\times \pts_2)$).
Reconsidering the analysis in (\ref{eq:lowQ}), we have
\begin{equation*}
|Q| = \sum_{i=1}^x |E_i|^2 > \sum_{i=1}^{n/2} d_i^2 = \Omega(n^3).
\end{equation*}
That is, if we perform the rest of the analysis as in the Elekes-Sharir framework, we only obtain
(at best) the weak lower bound $D(\pts')=\Omega(n)$.
The reason for this discrepancy is that the ``trivial'' distances $1,2,\ldots,n/2$
generate $\Omega(n^3)$ quadruples in $Q$, whereas the ``real'' inter-line distances
generate only a nearly quadratic number of quadruples. The discrepancy would have
been much lower if the analysis could have discarded the trivial distances, for example, by considering the \emph{bipartite} version of the problem, which only takes into account the distinct distances in $\pts_1\times\pts_2$; see \cite{SSS13} for a study of this variant.
The non-bipartite situation is another instance in which the values $|E_i|$ are (significantly more)
non-uniformly distributed.

\section{Distinct distances in rectangular lattices} \label{sec:rect}
In the second part of this note we consider the number of distinct distances that are
determined by an $n^{1-\alpha}\times n^{\alpha}$ integer lattice, for some
$0<\alpha\le 1/2$. We denote this number by $D_{\alpha}(n)$.

The case $\alpha=1/2$ is the case of the square $\sqrt{n}\times\sqrt{n}$
lattice, which determines $D_{1/2}(n)=\Theta(n/\sqrt{\log n})$ distinct distances, as reviewed
in Section \ref{sec:limitation}. Surprisingly, we show here a different estimate
 for $\alpha<1/2$.

%We first note that the upper bound holds for any $\alpha$, because any distance between
%a pair of points in the lattice is also realized as a distance from the lower-leftmost
%point of the lattice.

%The lower bound is trivial for $\alpha\le 1/3$; the easy argument is provided in the
%following lemma.
%
%
%
%%%%%%%%%%%%%%%%%%%%%%%%%%%%%%%%%%%%
%
%\begin{lemma}
%Let $R_\alpha(n)$ be an $n^{1-\alpha}\times n^{\alpha}$ integer lattice,
%for $0<\alpha\le 1/3$. Then the points of $\cal L$ determine $\Theta(n)$
%distinct distances.
%\end{lemma}
%
%%%%%%%%%%%%%%%%%%%%%%%%%%%%%%%%%%%%
%
%\begin{theProof}{\!\!}
%We establish the lower bound by considering only the following subset
%\[ D' = \{ \sqrt{a^2+b^2} \mid \tfrac12 n^{1-\alpha}\le a \le n^{1-\alpha}
%\text{ and } 0\le b\le n^{\alpha} \} \]
%of distinct distances, and notice that for any pair $(a,b)$ in the ranges
%specified for $D'$ we have $a^2+b^2<(a+1)^2$. Thus all the distances generated
%by these pairs are distinct, and so
%$D(R_\alpha(n))\ge|D'|\ge \frac12 n^{1-\alpha}\cdot n^\alpha = n/2$.
%\end{theProof}

%The analysis is more involved for the ``intermediate'' range
%$1/3 < \alpha < 1/2$, and is provided in the following main result of this section.
%
%
%
%As above, let $R_{\alpha}(n)$ denote the rectangular lattice
%$\{(i,j) \mid 0\le i\le n^{1-\alpha},\ 0\le j\le n^{\alpha}\}$,
%for some $1/3<\alpha \le 1/2$.
%
%We are interested in the cardinality of the set
%$$Q_{\alpha}(n)=\{i^2+j^2 \mid (i,j)\in R_{\alpha}(n)\} ,$$
%which is the number of distinct distances between the points of $R_\alpha(n)$.

%%%%%%%%%%%%%%%%%%%%%%%%%%

\begin{theorem}\label{Main}
For $\alpha<1/2$ we have $D_{\alpha}(n)=\Theta(n)$.
%$$|Q_{\alpha}(n)|=|R_{\alpha}(n)| - O(n^{2\alpha}\log^2n).$$
%In particular we have $|Q_{\alpha}(n)| = \Theta(n)$.
\end{theorem}

%%%%%%%%%%%%%%%%%%%%%%%%%%
\begin{theProof}{\!\!}
We consider the rectangular lattice
\[ R_{\alpha}(n) =\{(i,j)\ \mid \ 0\le i\le n^{1-\alpha},\ 0\le j\le n^{\alpha}\},\]
and its sublattice, 
\[R'_{\alpha}(n) =\{(i,j)\ \mid \  2n^{\alpha}\le i\le n^{1-\alpha},\ 0\le j\le n^{\alpha}\};\]
since $\alpha<1/2$, $R'_\alpha(n)\neq\emptyset$ for $n\ge n_0(\alpha)$, a suitable constant depending on $\alpha$. We also consider the functions
\begin{align*}
r(m) & = \bigl|\{(i,j)\in R'_{\alpha}(n)\ \mid \  i^2+j^2=m\}\bigr|, \\
d(m) & = \bigl|\{(i,j)\in R'_{\alpha}(n)\ \mid \ i^2-j^2=m\}\bigr|.
\end{align*}
Observe that the smallest (resp., largest) value of $m$ for which
$d(m)\ne 0$ is $3n^{2\alpha}$ (resp., $n^{2-2\alpha}$).

We have the identities
\begin{align}\label{re1}
\sum_m r(m) & = \sum_m d(m) , \\
\label{re2}
\sum_m r^2(m) & = \sum_m d^2(m) .
\end{align}
The identity (\ref{re1}) is trivial.
To see (\ref{re2}) we observe that the sum $\sum_m r^2(m)$
counts the number of ordered quadruples
$(i,j,i',j')$, for $(i,j),(i',j')\in R'_{\alpha}(n)$, such that $i^2+j^2=i'^2+j'^2$.
But this quantity also counts the number of those ordered quadruples
$(i,j,i',j')$, for $(i,j'),(i',j)\in R'_{\alpha}(n)$, such that $i^2-j'^2=i'^2-j^2$,
which is the value of the sum $\sum_m d^2(m)$. Putting (\ref{re1}) and
(\ref{re2}) together we have
\begin{equation}\label{re3}
\sum_m \binom{r(m)}2 = \sum_m \binom{d(m)}2.
\end{equation}
It is clear that $D_{\alpha}(n)\le |R_{\alpha}(n)| = n+O(n^{1-\alpha})$. In the rest of the proof we derive a matching lower bound for $D_{\alpha}(n)$.

Writing $M_k$ for the set of those $m$ with $r(m)=k$, we have
$\sum_k k|M_k|=|R'_{\alpha}(n)|$. On the other hand,
\begin{align*}
D_{\alpha}(n)& \ge  \sum_{k\ge 1}|M_k| \\
& = \sum_{k\ge 1} k|M_k|-\sum_{k\ge 1}(k-1)|M_k| \\
& = |R'_{\alpha}(n)|-\sum_{k\ge 2}(k-1)|M_k|.
\end{align*}

Thus
$D_{\alpha}(n)\ge n-O(n^{2\alpha})-\sum_{k\ge 2}(k-1)|M_k|$.
Using the inequality $k-1\le \binom k2$ and (\ref{re3}), we have
$$
\sum_{k\ge 2}(k-1)|M_k|\le \sum_{k\ge 2}\binom k2|M_k|=
\sum_m \binom{r(m)}2=\sum_m\binom{d(m)}2 .
$$
Theorem \ref{Main} is therefore a trivial consequence of the following proposition.
\end{theProof}

%%%%%%%%%%%%%%%%%%%%%%%%

\begin{proposition}\label{prop}
$$\sum_m\binom{d(m)}2 = O\left( n^{2\alpha}\log^2n\right).$$
\end{proposition}
%%%%%%%%%%%%%%%%%%%%%%%%
\begin{theProof}{\!\!}
We need the following easy lemma.
%%%%%%%%%%%%%%%%%%%%%%%%
\begin{lemma}\label{prod}
If $m$ can be written as the product of two integers in two different ways, say
$m=m_1m_2=m_3m_4$, then there exists a quadruple of positive integers $(s_1,s_2,s_3,s_4)$ satisfying
$$m_1=s_1s_2,\quad m_2=s_3s_4,\quad m_3=s_1s_3,\quad m_4=s_2s_4.$$
\end{lemma}
%%%%%%%%%%%%%%%%%%%%%%%%
\begin{theProof}{\!\!}
Since $m_1$ divides $m_3m_4$, we have $m_1=s_1s_2$ for some $s_1\mid m_3$ and
some $s_2\mid m_4$. Putting $s_3=m_3/s_1$ and $s_4=m_4/s_2$, we have
$m_2=s_3s_4$, $m_3=s_1s_3$, and $m_4=s_2s_4$.
\end{theProof}

We write
$$
\sum_m\binom{d(m)}2 = \sum_{1\le l\le n^{1-2\alpha}}
\sum_{m\in I_l}\binom{d(m)}2 ,
$$
where $I_l=[l^2n^{2\alpha},(l+1)^2n^{2\alpha})$.
We observe that the union of the intervals, namely
$[n^{2\alpha},(1+n^{1-2\alpha})^2n^{2\alpha})$, covers all the possible $m$ with $d(m)\ne 0$.

Now we estimate $\sum_{ m\in I_l}\binom{d(m)}2$ for a fixed $l$.
Let $a^2-b^2=c^2-d^2$ ($a>c$ and $b>d$) be a pair of distinct representations of some $m$,
which is counted in the above sum $\sum_{m\in I_l}\binom{d(m)}2$.
Since $m\in I_l$ we have
$$l^2n^{2\alpha}\le a^2-b^2< (l+1)^2n^{2\alpha}.$$
Thus,
$$
l^2n^{2\alpha}\le a^2< (l+1)^2n^{2\alpha}+b^2\le ((l+1)^2+1)n^{2\alpha} < (l+2)^2n^{2\alpha}.
$$
The same inequality holds for $c$, so we have
\begin{equation}\label{con1}
ln^{\alpha}\le a,c<(l+2)n^{\alpha}.
\end{equation}

Applying Lemma \ref{prod} to $(a-c)(a+c)=(b-d)(b+d)$, we obtain a quadruple
$(s_1,s_2,s_3,s_4)$ satisfying
\begin{align*}
& s_1s_2=a-c,\qquad s_3s_4=a+c,\\
& s_1s_3=b-d,\qquad s_2s_4=b+d.
\end{align*}

Using (\ref{con1}) and $0\le b,d\le n^{\alpha}$
we have the following inequalities:
\begin{align}
1\le s_1s_2,\ & s_1s_3,\ s_2s_4\le 2n^{\alpha},\nonumber\\
2ln^{\alpha}\le & s_3s_4< (2l+4)n^{\alpha}\label{des}.
\end{align}

It is clear from the above inequalities that $s_i\le 2n^{\alpha}$, for $i=1,\ldots,4$.
From $s_2s_4\le 2n^{\alpha},\ s_1s_3\le 2n^{\alpha}$, and $2ln^{\alpha}\le s_3s_4$,
we also deduce that
\begin{equation}\label{sa}
1\le s_2\le \frac{s_3}l \qquad\text{ and }\qquad
1\le s_1\le \frac{s_4}l.
\end{equation}
Choose $s_3$ between $1$ and $2n^{\alpha}$. Then choose $s_4$, according to (\ref{des}),
in the range $[\frac{2ln^{\alpha}}{s_3},\frac{(2l+4)n^{\alpha}}{s_3})$. Then choose
$s_1$ and $s_2$, according to (\ref{sa}), in
$\frac{s_3}l\cdot\frac{s_4}l\le \frac{(2l+4)n^{\alpha}}{l^2}$ ways.
The overall number of quadruples $(s_1,s_2,s_3,s_4)$ under consideration
is thus at most
$$\sum_{1\le s_3\le 2n^{\alpha}}\frac{4n^{\alpha}}{s_3}\cdot
\frac{(2l+4)n^{\alpha}}{l^2}=O\left (\frac{n^{2\alpha}\log n}l\right ).$$
Finally we have
$$\sum_m \binom{d(m)}2\le \sum_{1\le l\le n^{1-2\alpha}} \sum_{  m\in I_l}\binom{d(m)}2
= O\left(\sum_{l\le n^{1-2\alpha}}\frac{n^{2\alpha}\log n}l \right)
= O\left(n^{2\alpha}\log^2n\right).$$
\end{theProof}

\noindent{\bf Discussion.}
Theorem~\ref{Main} is closely related to a special case of a fairly deep conjecture
in number theory, stated as Conjecture 13~in Cilleruelo and Granville~\cite{CG07}.
This special case, given in \cite[Eq.~(5.1)]{CG07},
asserts that, for any integer $N$, and any fixed $\beta<1/2$,
$$
\bigl| \{ (a,b) \mid a^2+b^2=N,\; |b|<N^\beta \} \bigr| \le C_\beta ,
$$
where $C_\beta$ is a \emph{constant} that depends on $\beta$ (but not on $N$).
A simple geometric argument shows that this is true for $\beta\le 1/4$ but it is unknown for any $1/4<\beta<1/2$.
If that latter conjecture were true, a somewhat weaker version of
Theorem~\ref{Main} would follow.
Indeed, let $N$ be an integer that can be written as $i^2+j^2$, for
$\frac{1}{2}n^{1-\alpha} \le i\le n^{1-\alpha}$ and $j\le n^{\alpha}$. Then $N = \Theta(n^{2(1-\alpha)})$,
and $j = O(N^\beta)$, for $\beta=\alpha/(2(1-\alpha)) < 1/2$.

Conjecture 13 of \cite{CG07} would then imply that the number of pairs
$(i,j)$ as above is at most the constant $C_\beta$. In other words, each of the
$\Theta(n)$ distances in the portion of $R_\alpha(n)$ with $i\ge \frac{1}{2}n^{1-\alpha}$, interpreted as a distance from the origin
$(0,0)$, can be attained at most $C_\beta$ times.
Hence $D_\alpha(n) = \Theta(n)$, as asserted in Theorem \ref{Main}.

The general form of conjecture 13 \cite{CG07} asserts that the number of integer lattice points on an arc of length $N^{\beta}$ on the circle $a^2+b^2=N$ is bounded by some constant $C_{\beta}$, for any $\beta<1/2$. Cilleruelo and C\'ordoba \cite{CC} have proved this for $\beta<1/4$. See also Bourgain and Rudnick \cite{BR} for some consequences of this conjecture.

\vspace{3mm}

\noindent {\bf Acknowledgements.} The authors would like to thank Zeev Rudnick
for useful discussions on some of the number-theoretic issues. The third author
would like to thank Orit Raz and Joshua Zahl for several helpful discussions.

%%%%%%%%%%%%%%%%%%%%%%%%%%%%%%%%%%%%%%%%%%%%%%%%%%%%%%%%%%%%%%%%%%%%%%%%%%%%%%%%%

\end{document}